\documentclass[12pt]{amsart}
\usepackage{amsmath}
\usepackage{amssymb}
\usepackage{amsthm}

\newtheorem{theorem}{Theorem}[section]
\newtheorem{prop}[theorem]{Proposition}
\newtheorem{lem}[theorem]{Lemma}

\theoremstyle{definition}
\newtheorem{defn}[theorem]{Definition}
\theoremstyle{remark}
\newtheorem*{rem}{Remark}
\numberwithin{equation}{section}

\newcommand{\Gg}{\mathfrak{g}}    
\newcommand{\Gh}{\mathfrak{h}}
\newcommand{\Gq}{\mathfrak{q}}

\newcommand{\Gk}{\mathfrak{k}}

\newcommand{\Gz}{\mathfrak{z}}

\begin{document}

\title[Deformation Quantization of Poisson Structures]
{Deformation Quantization of\\ Certain Non-linear Poisson
Structures}

\author{Byung--Jay Kahng}
\date{}
\address{Department of Mathematics\\ University of California at Berkeley\\
Berkeley, CA 94704}
\email{bjkahng@math.berkeley.edu}

\begin{abstract}
 As a generalization of the linear Poisson bracket on the dual space
of a Lie algebra, we introduce certain non-linear Poisson brackets 
which are ``cocycle perturbations'' of the linear Poisson bracket. 
We show that these special Poisson brackets are equivalent to Poisson 
brackets of central extension type, which resemble the central extensions
of an ordinary Lie bracket via Lie algebra cocycles.  We are able to 
formulate (strict) deformation quantizations of these Poisson brackets 
by means of twisted group $C^*$--algebras.  We also indicate that
these deformation quantizations can be used to construct some specific 
non-compact quantum groups.
\end{abstract}
\maketitle

{\sc Introduction.} Let $M$ be a Poisson manifold.  Consider $C^{\infty}
(M)$, the commutative algebra under pointwise multiplication of smooth 
functions on $M$.  We attempt to deform the pointwise product of 
smooth functions into a noncommutative product, with respect to a 
parameter $\hbar$, such that the direction of the deformation is 
given by the Poisson bracket on $M$.  This problem of finding a 
{\em deformation quantization\/} of $M$ (\cite{Vy}, \cite{BFFL}) is 
actually a problem dating back to the early days of quantum mechanics
\cite{vN}, \cite{My}.

 We are particularly interested in the settings where the deformed 
product of functions is again a function---in contrast to much of 
the literature on the subject involving formal power series, or the 
so-called ``star products''.  In this direction, Rieffel has been
developing the notion of ``strict'' deformation quantization, in 
the $C^*$--algebra framework \cite{Rf1,Rf4}.  Here, in addition to 
the requirement that the deformed product of functions is again a 
function, the deformed algebra is required to have an involution 
and a $C^*$--norm.  By using the $C^*$--algebra framework, one 
gains the advantage of being able to keep the topological and 
geometric aspects of the given manifold while we perform the 
quantization.

 Let $\Gh$ be a (finite dimensional) Lie algebra.  It is well-known
\cite{Wn1} that the Lie algebra structure on $\Gh$ defines a natural
Poisson bracket on the dual vector space $\Gh^*$ of $\Gh$, which is
called a {\em linear Poisson bracket\/}.  This Poisson bracket is 
also called the ``Lie--Poisson bracket'', to emphasize the fact that
it actually was already known to Lie.  In \cite{Rf3} Rieffel showed 
that given the linear Poisson bracket on $\Gh^*$, a deformation 
quantization is provided by the convolution algebra structure on the 
simply connected Lie group $H$ corresponding to $\Gh$.  In particular,
when $\Gh$ is a nilpotent Lie algebra, this is shown to be a strict 
deformation quantization.

 In this paper, we wish to generalize the above situation and to
include twisted group convolution algebras into the framework of
deformation quantization.  We first define a class of Poisson 
brackets on the dual vector space of a Lie algebra, which contains 
the linear Poisson bracket as a special case.  These Poisson brackets 
can actually be realized as ``central extensions'' of the linear 
Poisson bracket.  We then show that twisted group convolution 
algebras provide deformation quantizations of these Poisson brackets.
We obtain strict deformation quantizations when the Lie algebra is 
nilpotent. 

 In addition to its interest as a generalization of the deformation
quantization problem into the non-linear situation, this result has
a nice application to quantum groups.  Quantum groups \cite{Dr}, 
\cite{CP} are usually obtained by suitably ``deforming'' ordinary 
Lie groups, and as suggested by Drinfeld \cite{Dr}, we expect to 
obtain quantum groups by deformation quantization of the so-called 
Poisson--Lie groups \cite{LW}.  In some cases, the compatible Poisson
brackets on the Poisson--Lie groups are shown to be of our special 
type, in which case we can apply the result of this paper to obtain
(strict) deformation quantizations of them.  

 This enables us to construct some specific non-compact quantum 
groups.  Not only have we actually been able to show \cite{BJK} that 
some of the earlier known examples of non-compact quantum groups 
\cite{Rf5}, \cite{SZ}, \cite{VD}, \cite{Ld} are obtained in this way,
but we also obtain a new class of non-compact quantum groups \cite{BJK2}.
Although the method of construction may seem rather naive, our new 
example is shown to satisfy some interesting properties, including the
``quasitriangular'' property.  We will discuss our construction of 
quantum groups in a separate paper.

 This paper is organized as follows.  In the first section, we 
review the definitions of Poisson brackets and the formulation 
of (strict) deformation quantization.  We also include a 
discussion on twisted group algebras, which are the main 
objects of our study.  In the second section, we define
our special class of non-linear Poisson brackets, as motivated by
the central extension of ordinary Lie brackets.  We show in the 
third section that certain twisted group $C^*$--algebras provide 
strict deformation quantizations of these special Poisson 
brackets.  We use some non-trivial results on twisted group 
$C^*$--algebras obtained by Packer and Raeburn \cite{PR1, PR2}. 
We restrict our study to the strict deformation quantization 
case, but some indications for generalization are also briefly 
mentioned.  
  
 The essential part of this article is from the author's
Ph.D. thesis at U.C. Berkeley.  I would like to express here
my deep gratitude to Professor Marc Rieffel, without whose
constant encouragement, show of interest and numerous 
suggestions, this work would not have been made possible.

\section{Preliminaries}

 Let $M$ be a $C^\infty$ manifold, and let $C^\infty(M)$ be the 
algebra of complex-valued $C^\infty$ functions on $M$.  It is a 
commutative algebra under pointwise multiplication, and is equipped
with an involution given by complex conjugation. 

\begin{defn}\label{2.2.1}
By a {\it Poisson bracket\/} on $M$, we mean a skew, bilinear map 
$\{\ ,\ \}:C^\infty(M)\times C^\infty(M)\rightarrow C^\infty(M)$ such 
that the following holds:

\begin{itemize}
\item $\{\ ,\ \}$ defines a Lie algebra structure on $C^\infty(M)$.

 (i.\,e.\ the bracket satisfies the Jacobi identity.)
\item (Leibniz rule): $\{f,gh\}=\{f,g\}h+g\{f,h\}$, for $f,g,h\in
C^\infty(M)$.
\end{itemize}
We also require that the Poisson bracket be real, in the sense that 
$\{f^*,g^*\}=\{f,g\}^*$.  A manifold $M$ equipped with such a bracket
is called a {\it Poisson manifold\/}, and $C^\infty(M)$ is a {\it 
Poisson (*--)algebra\/}.
\end{defn}

 The deformation quantization will take place in $C^\infty(M)$ (or to
allow non-compact $M$, in $C^\infty_\infty(M)$, which is the space of
$C^\infty$ functions vanishing at infinity).  The Poisson bracket on $M$
gives  the direction of the deformation.  Let us formulate the  following 
definition for deformation quantization, which is the  ``strict'' 
deformation quantization proposed by Rieffel \cite{Rf1}. Depending on 
the situations, we may also be interested in some other subalgebras of
$C^\infty(M)$, on which the deformation takes place.   So the definition 
is formulated at the level of an arbitrary (dense) ${}^*$--subalgebra
$\mathcal A$, for example the algebra of Schwartz  functions. 

\begin{defn}\label{2.2.2} 
Let $M$ be a Poisson manifold as above. Let $\mathcal A$ be a dense 
$*$--subalgebra (with respect to the $C^*$--norm $\|\ \|_{\infty}$) 
of $C_{\infty}(M)$, on which $\{\ ,\ \}$ is defined with values in 
$\mathcal A$.  By a {\it strict deformation quantization\/} of 
$\mathcal A$ in the direction of $\{\ ,\ \}$, we mean an open 
interval $I$ of real numbers containing $0$, together with a triple 
($\times_\hbar$,${}^{*_\hbar}$, $\|\ \|_\hbar$) for each $\hbar\in I$,
of an associative product, an involution, and a $C^*$--norm (for 
$\times_\hbar$ and ${}^{*_\hbar}$) on $\mathcal A$, such that

\begin{enumerate} 
\item For $\hbar=0$, the operations $\times_\hbar$,${}^{*_\hbar}$,
$\|\ \|_\hbar$ are the original pointwise product, involution (complex
conjugation), and $C^*$--norm (i.\,e.\ sup-norm $\|\ \|_\infty$) on 
$\mathcal A$, respectively. 
\item The completed $C^*$-algebras $A_\hbar$ form a ``continuous 
field'' of $C^*$--algebras (In particular, the map $\hbar \mapsto 
\|f\|_\hbar$ is continuous for any $f\in {\mathcal A}$.).
\item For any $f,g\in {\mathcal A}$,
$$ 
\left\|\frac{f\times_\hbar g - g\times_\hbar f}{i\hbar}-\{f,g\}
\right\|_\hbar\to 0 
$$
as $\hbar \to 0$.
\end{enumerate} 
\end{defn}

 Our main result below will show that certain ``twisted group 
$C^*$--algebras'' provide strict deformation quantizations.  Let us
briefly discuss about these algebras, mainly to set up our notation.
We give here only those results that are needed in later sections.
For more discussion on the subject, we refer the reader to the 
articles by Zeller-Meier \cite{ZM} (when the group is discrete) and
by Busby and Smith \cite{BuS}.  There are also many other recent 
articles on these algebras including \cite{PR1,PR2,PR3}, which contain
considerably deeper results.  To avoid technical pathology, we assume
that the groups we consider are discrete or second countable locally 
compact (e.\,g.\ Lie groups), and the $C^*$--algebras are always
separable.

\begin{defn}\label{2.3.8}(\cite{Rf2})
Let $G$ be a locally compact group with left Haar measure $dx$ and 
modular function ${\Delta}_G$.  Let $G$ act on a $C^*$--algebra $A$ 
and let us denote this action by $\alpha:G\to \operatorname{Aut}(A)$. 
We assume that $\alpha$ is strongly continuous.  Also let $UZM(A)$ 
be the group of unitary elements in the center of the multiplier 
algebra, $M(A)$, of $A$.  By a {\em continuous field\/} over $N$ of 
{\em $\alpha$--cocycles\/} of $G$, where $N$ is a locally compact 
space, we will mean a function $\sigma$ on $G\times G\times N$ with 
values in $UZM(A)$ such that
\begin{itemize}
\item If we fix $r\in N$, then $\sigma$ is a normalized 
$\alpha$--cocycle on $G$. That is,
$$
\bigl({\alpha}_x\sigma(y,z;r)\bigr)\sigma(x,yz;r)=
\sigma(x,y;r)\sigma(xy,z;r)
$$
and
$$
\sigma(x,e;r)=\sigma(e,x;r)=1
$$
for $x,y,z\in G$, where $e$ is the identity element of $G$. 
\item If we fix $x,y\in G$, then $\sigma$ is continuous on $N$.
\item For any $f\in C_{\infty}(N,A)$ the function
$$
(x,y)\mapsto f(\cdot)\sigma(x,y;\cdot)
$$
from $G\times G$ to $C_{\infty}(N,A)$ is Bochner measurable.
\end{itemize}
\end{defn}

 For convenience, let us denote by $\sigma^r,r\in N$ the ordinary group
cocycle on $G$ defined by $\sigma^r(x,y)=\sigma(x,y;r)$.  Corresponding 
to the continuous field of cocycles $\sigma:r\mapsto\sigma^r$, we can 
define \cite{Rf2} the twisted convolution and involution on $L^1\bigl
(G,C_{\infty}(N,A)\bigr)$.  We define, for $\phi,\psi\in L^1\bigl
(G,C_{\infty}(N,A)\bigr)$:
$$
(\phi*_{\sigma}\psi)(x;r)=\int_G \phi(y;r){\alpha}_y\bigl(\psi(y^{-1}x;r)
\bigr){\sigma}^{r}(y,y^{-1}x)\,dy
$$
and
$$
{\phi}^{*}(x;r)={\alpha}_x\bigl(\phi(x^{-1};r)^{*}\bigr){\sigma}^r
(x,x^{-1})^{*}{\Delta}_G(x^{-1}).
$$
We thus obtain a Banach ${}^*$--algebra.  Let us denote this algebra
by $L^1(G,N,A,\sigma)$, with the group action $\alpha$ to be understood.
We also define $C^*(G,N,A,\sigma)$, the enveloping $C^*$--algebra of 
$L^1(G,N,A,\sigma)$.  There are also the notions of induced 
representations and regular representations \cite{BuS}, \cite{ZM}, 
\cite{Rf2}.  So we may as well define the reduced $C^*$--algebra $C^*_r
(G,N,A,\sigma)$.  All these are more or less straightforward.

 Compare this definition with the definition given in \cite{BuS}, where
the twisted group convolution algebra has been formulated via a single
cocycle.  Nevertheless, the present definition is no different from the
usual one, since we can regard $\sigma$ also as a single cocycle, taking
values in $UZM\bigl(C_{\infty}(N,A)\bigr)$.  The present formulation is
useful when we study the continuity problem of the fields of 
$C^*$--algebras consisting of twisted group $C^*$--algebras, by varying
the cocycles.  Recall that the continuous field property is essential 
in the definition of the strict deformation quantization (Definition 
\ref{2.2.2}).  

 Using the universal property of full $C^*$--algebras, and also taking 
advantage of the property of the reduced $C^*$--algebras that one is 
able to work with their specific representations, Rieffel in \cite{Rf2} 
gave an answer to the problem of the continuity of the field of 
$C^*$--algebras $\bigl\{C^*(G,A,{\sigma}^{r})\bigr\}_{r\in N}$, as 
follows.  Here $C^*(G,A,\sigma^r)$ is the twisted group $C^*$--algebra 
in the usual sense of \cite{BuS}.  

\begin{theorem}\label{2.3.9}
Let $G,A,\alpha$ be understood as above.  Let $\sigma$ be a 
continuous field over $N$ of $\alpha$--cocycles on $G$.  Then
\begin{itemize}
\item The field $\bigl\{C^*(G,A,{\sigma}^{r})\bigr\}$ over $N$ is 
upper semi-continuous.
\item $\bigl\{C^*_r(G,A,{\sigma}^{r})\bigr\}$ over $N$ is lower 
semi-continuous.
\item Thus, if each $(G,A,\alpha,{\sigma}^{r})$ satisfies the 
``amenability condition'', i.\,e.\ $C^*_r(G,A,\alpha,\sigma^r)=
C^*(G,A,\alpha,\sigma^r)$, then it follows that the field of 
$C^*$--algebras $\bigl\{C^*(G,A,{\sigma}^{r})\bigr\}_{r\in N}$ is 
continuous.
\end{itemize}
\end{theorem}

 For the proof of the theorem and the related questions, we will 
refer the reader to \cite{Rf2}, and the references therein.  Note 
that by replacing $A$ with $C_{\infty}(N,A)$ and by introducing a
new base space, we may even consider a continuous field of the 
twisted group $C^*$--algebras given by the cocycles of continuous
field type (Definition \ref{2.3.8}).   

 The $C^*$--algebra $C^*(G,N,A,\sigma)$ may be regarded as a 
$C^*$--algebra of ``cross sections'' of the continuous field $\bigl
\{C^*(G,A,{\sigma}^{r})\bigr\}_{r\in N}$.  It is called the {\em 
$C^*$--algebra of sections of a $C^*$--bundle\/} by Packer and Raeburn
\cite{PR2} (Compare this terminology with Fell's notion of 
``$C^*$--algebraic bundles'' \cite{FD}, which is considerably more 
general than is needed for our present purposes.).  Actually in 
\cite{PR2}, the continuity problem of twisted group $C^*$--algebras 
allowing both the cocycle and the action to vary continuously has been
studied in terms of the aforementioned notion of section $C^*$--algebra
of a $C^*$--bundle.  Taking a related viewpoint, Blanchard in \cite{Bl}
has recently developed a framework for a general continuous field of 
$C^*$--algebras in terms of ``$C_{\infty}(X)$--algebras'', where $X$ 
in our case is the locally compact base space $N$.  A $C_{\infty}
(X)$--algebra is a certain $C^*$--algebra having a $C(X)$ module 
structure.  See \cite{Bl}.    

 We conclude this section by quoting (without proof) a couple 
of deep theorems of Packer and Raeburn \cite{PR1,PR2} on the 
structure of twisted group $C^*$--algebras.  We tried to keep Packer 
and Raeburn's notation and terminology.  Although some of them are 
different from our notation, they are clear enough to understand.  
For example, $A\times_{\alpha,u}G$ denotes the twisted group 
$C^*$--algebra (or ``twisted crossed product'') $C^*(G,A,\alpha,u)$. 
All this and more can be found in \cite{PR1,PR2}.  These theorems 
will be used later in the proof of our Theorem \ref{2.4.8}, which is
our main result.

\begin{theorem}\label{b1}(\cite{PR1}) 
(Decomposition of twisted crossed products)
Suppose that $(A,G,\alpha,u)$ is a separable twisted dynamical 
system and $N$ is a closed normal subgroup of $G$.  There exists a 
canonically determined twisted action $(\beta,v)$ of $G/N$ on 
$A\times_{\alpha,u}N$ such that:
$$
A\times_{\alpha,u}G\cong(A\times_{\alpha,u}N)\times_{\beta,v}G/N.
$$
\end{theorem}

 The next theorem is about the continuity of a field of twisted group
$C^*$--algebras.  Compare this with Theorem \ref{2.3.9}, where we 
considered the continuity problem only when the twisting cocycle is
varying.  Meanwhile, note in the theorem that $G$ is assumed 
to be amenable (So by the ``stabilization trick'' of Packer and
Raeburn \cite{PR1}, the amenability condition always holds for any
quadruple $(G,A,\alpha,u)$.).   

\begin{theorem}\label{b2}(\cite{PR2})
Suppose $A$ is the $C^*$--algebra of sections of a separable 
$C^*$--bundle over a locally compact space $X$, and $(\alpha,u)$ is a 
twisted action of an amenable locally compact group $G$ on $A$ such that
each ideal $I_x=\{a\in A:a(x)=0\}$ is invariant.  Then for each $x\in 
X$, there is a natural twisted action $\bigl(\alpha(x),u(x)\bigr)$ on
the quotient $A/{I_x}$, and $A\times_{\alpha,u}G$ is the $C^*$--algebra 
of sections of a $C^*$--bundle over $X$ with fibers isomorphic to $\left
(A/{I_x}\right)\times_{\alpha(x),u(x)}G$.
\end{theorem}

\section{The non-linear Poisson bracket}

 Let us begin by trying to characterize the special Poisson brackets 
that will allow twisted group algebras to be deformation quantizations 
of them.  Recall that the twisting of the convolution algebra structure 
in a twisted group algebra is given by group (2--)cocycles.  Meanwhile 
any group cocycle for a locally compact group $G$ having values in an 
abelian group $N$ can be canonically associated with a central extension
of $G$ by $N$, and actually all central extensions are essentially 
obtained in this way \cite{FD}.

 Since it is known \cite{Rf3} that ordinary group convolution algebras 
can be regarded as deformation quantizations of linear Poisson 
brackets, the above observations suggest that twisted group algebras 
will provide deformation quantizations of certain Poisson brackets 
which are, in a loose sense, ``central extensions'' of linear Poisson 
brackets.  Although we have to make clear what we mean by this last 
statement, this is the main motivation behind the definition of our 
special type of Poisson bracket formulated below.

 Let $\Gh$ be a (finite--dimensional) Lie algebra and let us denote by
$\Gg=\Gh^*$ its dual vector space.  As usual, we will denote the dual 
pairing between $\Gh$ and $\Gg$ by $\langle\ ,\ \rangle$.  We will 
have to require later that $\Gh$ is a nilpotent or an exponential 
solvable Lie algebra because of some technical reasons to be discussed
below, but for the time being we allow $\Gh$ to be a general Lie algebra.  
Recall \cite{Wn1} that we define the linear Poisson bracket 
on the dual vector space $\Gg=\Gh^*$ by
\begin{equation}\label{(2.2)}
\{\phi,\psi\}_{\text {lin}}(\mu)=\bigl\langle[d\phi(\mu),d\psi(\mu)],
\mu\bigr\rangle
\end{equation}
where $\phi,\psi\in C^{\infty}(\Gg)$ and $\mu\in\Gg$.  Here $d\phi(\mu)$
and $d\psi(\mu)$ has been naturally realized as elements in $\Gh$.

 We wish to define a generalization of this Poisson bracket by allowing
a suitable ``perturbation'' of the right-hand side of equation 
\eqref{(2.2)}.  This will be done via a certain Lie algebra 2--cocycle 
on $\Gh$, denoted by $\Omega$, having values in $C^{\infty}(\Gg)$.  That
is, we will consider the Poisson brackets of the form:
\begin{equation}\label{(2.3)}
\{\phi,\psi\}(\mu)=\bigl\langle[d\phi(\mu),d\psi(\mu)],\mu\bigr\rangle
+\Omega\bigl(d\phi(\mu),d\psi(\mu);\mu\bigr).
\end{equation}
As above, we regard $d\phi(\mu)$ and $d\psi(\mu)$ as elements in $\Gh$.
  
 Compare equation \eqref{(2.3)} with the definition of the linear 
Poisson bracket.  In the linear Poisson bracket case, the Lie bracket 
takes values in $\Gh$, the elements of which can be regarded as (linear)
functions contained in $C^{\infty}(\Gg)$, via the dual pairing.  That is,
the right-hand side of equation \eqref{(2.2)} can be viewed as the 
evaluation at $\mu\in\Gg$ of a $C^{\infty}$--function, $[X,Y]\in\Gh
\subseteq C^{\infty}(\Gg)$, where $X=d\phi(\mu)$ and $Y=d\psi(\mu)$. 
In the ``perturbed'' case, the right-hand side of equation \eqref{(2.3)}
may be viewed as the evaluation at $\mu\in\Gg$ of a 
$C^{\infty}$--function, $[X,Y]+\Omega(X,Y)\in C^{\infty}(\Gg)$, where 
$X=d\phi(\mu)$ and $Y=d\psi(\mu)$.  So to make sense of the Poisson 
brackets of the type given by equation \eqref{(2.3)}, we will first study
the ``perturbation'' of the Lie bracket on $\Gh$ by a cocycle.  Later, 
we will find some additional conditions for the cocycle $\Omega$ such 
that the equation \eqref{(2.3)} indeed gives a well-defined Poisson 
bracket on $\Gg$.

 Let $V$ be a $U(\Gh)$--module, possibly infinite dimensional.  Consider
a {\em 2--cocycle\/} $\Omega$ for $\Gh$ having values in $V$.  It is a
skew-symmetric, bilinear map from $\Gh\times\Gh$ into $V$ such that 
$d\Omega=0$ (For more discussion on cohomology of Lie algebras, see the 
standard textbooks on the subject \cite[\S5]{CE}, \cite{Kn}.).  When $V$
is further viewed as an abelian Lie algebra, the space $\Gh\oplus V$ can 
be given a Lie algebra structure \cite{B}, \cite{Kn} which becomes a 
{\em central extension\/} Lie algebra of $\Gh$ by $V$:
$$
\bigl[(X,v),(Y,w)\bigr]_{\Gh\oplus V}=\bigl([X,Y],X\cdot w-Y\cdot v+
\Omega(X,Y)\bigr)
$$
for $X,Y\in\Gh$ and $v,w\in V$.  Here the dot denotes the module action.
In particular, when $V$ is assumed to be a trivial $U(\Gh)$--module, we
have:
\begin{equation}\label{(2.5)}
\bigl[(X,v),(Y,w)\bigr]_{\Gh\oplus V}=\bigl([X,Y],\Omega(X,Y)\bigr).
\end{equation}

 Let us slightly modify this ``central extension'' picture as follows,
so that we are able to consider a Lie bracket on $\Gh+V$, where we now
allow $\Gh\cap V\ne{0}$ in general.  Clearly, $\Gh\cap V$ is a subspace
of $\Gh$.  However, since $\Gh$ is already equipped with its given Lie
bracket and since $V$ will be assumed to be an abelian Lie algebra, it
is only reasonable to consider the case in which $\Gh\cap V$ is an 
abelian subalgebra of $\Gh$.  For simplicity, we will further assume 
that $\Gh\cap V$ is a central subalgebra of $\Gh$, which means that $V$
is a trivial $U(\Gh)$--module.  Let us denote this central subalgebra
by $\Gz$.  Without loss of generality, we may assume that $\Gz$ is the 
center of $\Gh$.  In this case, we just replace $V$ by an extended 
abelian Lie algebra, still denoted by $V$, satisfying $\Gh\cap V=\Gz$. 

 Let us look for a trivial $U(\Gh)$--module $V$, which we will view as 
an abelian Lie algebra, such that $\Gh\cap V=\Gz$ is the center of $\Gh$.
Since we eventually want to define a $V$--valued cocycle for $\Gh$, 
from which we construct a bracket operation on $C^{\infty}(\Gg)$, we also
require that $V$ is contained in $C^{\infty}(\Gg)$.  So let us consider
the subspace $\Gq=\Gz^{\bot}$ of $\Gg$, and choose as our $V$ the 
following:
$$
V=C^{\infty}(\Gg/\Gq)\subseteq C^{\infty}(\Gg).
$$
Here the functions in $V=C^{\infty}(\Gg/\Gq)$ have been realized as
functions in $C^{\infty}(\Gg)$, by the ``pull-back'' using the natural
projection $p:\Gg\to\Gg/\Gq$.  

 Since any $X\in\Gh$ can be regarded as a linear function contained in
$C^{\infty}(\Gg)$ via the dual pairing, we can see easily that $X\in\Gh
\cap V\subseteq C^{\infty}(\Gg)$ if and only if $\langle X,\mu+\nu\rangle
=\langle X,\mu\rangle$ for all $\mu\in\Gg,\nu\in\Gq$.  It follows 
immediately that $\Gh\cap V=\Gz$.  On the other hand, consider the 
representation $\operatorname{ad}^*$.  For any $X\in\Gh$ and any $\mu
\in\Gg$, we have $\operatorname{ad}^*_{-X}(\mu)=\nu\in\Gq$, since for any
$Y\in\Gz$, we have $\langle Y,\nu\rangle=\bigl\langle Y,\operatorname{ad}
^*_{-X}(\mu)\bigr\rangle=\bigl\langle[X,Y],\mu\bigr\rangle=0$.  It follows
that 
$$
\operatorname{ad}^*_X(f)(\mu)=f\bigl(\operatorname{ad}^*_{-X}(\mu)
\bigr)=f(\nu)=0,
$$
for any $f\in V=C^{\infty}(\Gg/\Gq)$.  By the natural extension of 
$\operatorname{ad}^*$ to $U(\Gh)$, we can give $V$ the trivial 
$U(\Gh)$--module structure.

\begin{rem}
When $\Gh$ has a trivial center, the space $V$ will be just $\{0\}$.
To avoid this problem, we could have considered $C^{\infty}(\Gg)^H$,
the space of $\operatorname{Ad}^*H$--invariant $C^{\infty}$ functions 
on $\Gg$.  It is always nonempty (It contains the so-called Casimir 
elements \cite{V}.).  It also satisfies $\Gh\cap C^{\infty}(\Gg)^H=\Gz$
and can be given the trivial $U(\Gh)$--module structure.  However, it 
does not satisfy the following property (Lemma \ref{2.4.2}), which we 
need later when we define our Poisson bracket.  For this reason, we 
choose our $V$ as it is defined above.  At least for nilpotent $\Gh$, 
which is the case we are going to study most of the time, this is less
of a problem since $\Gh$ has a non-trivial center.  
\end{rem}

\begin{lem}\label{2.4.2}
Let $V\subseteq C^{\infty}(\Gg)$ be defined as above.  Then for any
function $\chi\in V$ and for any $\mu\in\Gg$, we have: $d\chi(\mu)\in
\Gz$.
\end{lem}

\begin{proof}
Since $V\subseteq C^{\infty}(\Gg)$, it follows that $d\chi(\mu)\in\Gh$.
Recall that $d\chi(\mu)$ defines a linear functional on $\Gg=\Gh^*$ by
$$
\bigl\langle d\chi(\mu),\nu\bigr\rangle=\left.\left(\frac{d}{dt}\right)
\right|_{t=0}\chi(\mu+t\nu).
$$
To see if $d\chi$ is contained in $\Gz$, suppose $\nu\in\Gq=\Gz^{\bot}$.
Since $\chi\in V=C^{\infty}(\Gg/\Gq)$, we know that $\chi(\mu+\nu)=
\chi(\mu)$, for all $\nu\in\Gq$.  Therefore, the above expression 
becomes:
$$
\bigl\langle d\chi(\mu),\nu\bigr\rangle=\left.\left(\frac{d}{dt}\right)
\right|_{t=0}\chi(\mu)=0.
$$
Since $\nu\in\Gq$ is arbitrary, we thus have: $d\chi(\mu)\in\Gz$.
\end{proof}

 Let us now turn to the discussion of defining a (perturbed) bracket
operation on $\Gh+V$, which will enable us to formulate our special
Poisson bracket on $C^{\infty}(\Gg)$.  Since $\Gh$ and $V$ are 
subspaces of $\Gh+V$, there exists a (linear) surjective map, $\Gh\to
(\Gh+V)/V$, whose kernel is $\Gh\cap V=\Gz$.  We thus obtain a vector 
space isomorphism, in a canonical way, between $(\Gh+V)/V$ and $\Gh/\Gz$. 
The map from $\Gh+V$ onto $\Gh/\Gz$ is a canonical one, which extends 
the canonical projection of $\Gh$ onto $\Gh/\Gz$.  Therefore, it is 
reasonable to consider a cocycle for $\Gh/\Gz$ having values in $V$ 
(viewed as a trivial $U(\Gh/\Gz)$--module) and use it to define a 
bracket operation on $\Gh+V$.  Let $\Omega$ be such a cocycle for 
$\Gh/\Gz$. 

\begin{rem}
Note that in this setting, the cocycle $\Omega$ can naturally be 
identified with a cocycle $\tilde{\Omega}$ for $\Gh$ having values 
in $V$ (considered as a trivial $U(\Gh)$--module), satisfying the
following ``centrality condition'':
\begin{equation}\label{(central)}
\tilde{\Omega}(Z,Y)=\tilde{\Omega}(Y,Z)=0,
\end{equation}
for any $Z\in\Gz$ and any $Y\in\Gh$.  In fact, we may define $\tilde
{\Omega}$ as $\tilde{\Omega}(X,Y)=\Omega(\dot{X},\dot{Y})$, where
$\dot{X}$ denotes the image in $\Gh/\Gz$ of $X$ under the canonical
projection.  For this reason, we will from time to time use the same
notation, $\Omega$, to denote both $\Omega$ and $\tilde{\Omega}$.
\end{rem}

 By viewing $\Omega$ as a cocycle for $\Gh$, we can define, as in 
equation \eqref{(2.5)}, a Lie bracket on $\Gh\oplus V$:
$$
\bigl[(X,v),(Y,w)\bigr]_{\Gh\oplus V}=\bigl([X,Y],\Omega(X,Y)\bigr).
$$
To define a bracket operation on $\Gh+V$, consider the natural 
surjective map from $\Gh\oplus V$ onto $\Gh+V$, whose kernel is:
$$
\delta=\{(Z,-Z):Z\in\Gz\}\subseteq\Gh\oplus V.
$$
Since $\delta$ is clearly central with respect to the Lie bracket $[\ ,\
]_{\Gh\oplus V}$ given above, it is an ideal.  Therefore, $\Gh+V=(\Gh
\oplus V)/{\delta}$ is a Lie algebra.  The Lie bracket on it is given
by:
\begin{equation}\label{(2.6)}
[X+v,Y+w]_{\Gh+V}=[X,Y]+\Omega(X,Y),\qquad X,Y\in\Gh,\quad v,w\in V
\end{equation}
which is the given Lie bracket on $\Gh$ plus a cocycle term.  In this
sense, equation \eqref{(2.6)} may be considered as a ``perturbed Lie 
bracket'' of the given Lie bracket on $\Gh$.  Compare this with equation 
\eqref{(2.3)}, where the linear Poisson bracket (given by the Lie bracket
on $\Gh$) is ``perturbed'' by a certain cocycle $\Omega$. 

 Using the observation given above as motivation, let us define more 
rigorously our Poisson bracket on $\Gg=\Gh^*$.  This is, in fact, a 
``cocycle perturbation'' of $\{\ ,\ \}_{\text {lin}}$ on $\Gg$.

\begin{theorem}\label{2.4.3}
Let $\Gh$ be a Lie algebra with center $\Gz$ and let $\Gg=\Gh^*$ be the 
dual vector space of $\Gh$.  Consider the vector space $V=C^{\infty}
(\Gg/\Gq)\subseteq C^{\infty}(\Gg)$ as above, where $\Gq=\Gz^{\bot}$.
Let us give $V$ the trivial $U(\Gh)$--module structure.  Let $\Omega$ be
a Lie algebra 2--cocycle for $\Gh$ having values in $V$, satisfying the
centrality condition.  That is, $\Omega$ is a skew-symmetric, bilinear 
map from $\Gh\times\Gh$ into $V$ such that:
$$
\Omega\bigl(X,[Y,Z]\bigr)+\Omega\bigl(Y,[Z,X]\bigr)+\Omega\bigl(Z,[X,
Y]\bigr)=0,\qquad X,Y,Z\in\Gh
$$
satisfying: $\Omega(Z,Y)=\Omega(Y,Z)=0$ for $Z\in\Gz$ and any $Y\in\Gh$.
Then the bracket operation $\{\ ,\ \}_{\Omega}:C^{\infty}(\Gg)\times 
C^{\infty}(\Gg)\to C^{\infty}(\Gg)$ defined by
$$
\{\phi,\psi\}_{\Omega}(\mu)=\bigl\langle[d\phi(\mu),d\psi(\mu)],\mu
\bigr\rangle+\Omega\bigl(d\phi(\mu),d\psi(\mu);\mu\bigr)
$$
is a Poisson bracket on $\Gg$.  
\end{theorem}

\begin{rem}
If we denote $d\phi(\mu)$ and $d\psi(\mu)$ by $X$ and $Y$, as elements 
in $\Gh$, the right-hand side of the definition of the Poisson bracket
may be viewed as the evaluation at $\mu\in\Gg$ of a 
$C^{\infty}$--function, $[X,Y]+\Omega(X,Y)\in\Gh+V\in C^{\infty}(\Gg)$.
Note that this expression is just the Lie bracket on $\Gh+V$ defined 
earlier by equation \eqref{(2.6)}. 
\end{rem}

\begin{proof}
Since $d\phi(\mu)$ and $d\psi(\mu)$ can be naturally viewed as elements
in $\Gh$, it is easy to see that $\{\ ,\ \}_{\Omega}$ is indeed a map 
from $C^{\infty}(\Gg)\times C^{\infty}(\Gg)$ into $C^{\infty}(\Gg)$.  
The skew-symmetry and bilinearity are clear.

To verify the Jacobi identity, consider the functions $\phi_1,
\phi_2,\phi_3$ in $C^{\infty}(\Gg)$.  We may write $\{\phi_2,\phi_3\}
_{\Omega}$ as:
$$
\{\phi_2,\phi_3\}_{\Omega}(\mu)=\{\phi_2,\phi_3\}_{\text {lin}}(\mu)
+\chi(\mu),
$$
where $\chi$ is a function in $V$.  We therefore have:
$$
d\bigl(\{\phi_2,\phi_3\}_{\Omega}\bigr)(\mu)=\bigl[d\phi_2(\mu),
d\phi_3(\mu)\bigr]+d\chi(\mu).
$$
The first term in the right hand side is the differential of the linear
Poisson bracket, which is rather well known \cite{Wn1}.  Moreover, since 
$\chi\in V$, it follows from Lemma \ref{2.4.2} that $d\chi(\mu)\in\Gz$, 
which is ``central'' with respect to both $[\ ,\ ]$ and $\Omega$.  We 
thus have:
\begin{align}
\bigl\{\phi_1,&\{\phi_2,\phi_3\}_{\Omega}\bigr\}_{\Omega}(\mu) \notag \\
&=\bigl\langle\bigl[d\phi_1(\mu),d(\{\phi_2,\phi_3\}_{\Omega})(\mu)
\bigr],\mu\bigr\rangle+\Omega\bigl(d\phi_1(\mu),d(\{\phi_2,\phi_3\}
_{\Omega})(\mu);\mu\bigr)  \notag \\
&=\bigl\langle\bigl[d\phi_1(\mu),[d\phi_2(\mu),d\phi_3(\mu)]\bigr],\mu
\bigr\rangle+\Omega\bigl(d\phi_1(\mu),[d\phi_2(\mu),d\phi_3(\mu)];\mu
\bigr), \notag
\end{align}
and similarly for $\bigl\{\phi_2,\{\phi_3,\phi_1\}_{\Omega}\bigr\}
_{\Omega}$ and $\bigl\{\phi_3,\{\phi_1,\phi_2\}_{\Omega}\bigr\}
_{\Omega}$.  So the Jacobi identity for $\{\ ,\ \}_{\Omega}$ follows 
from that of the Lie bracket $[\ ,\ ]$ and the cocycle identity for 
$\Omega$.  That is,
$$
\bigl\{\phi_1,\{\phi_2,\phi_3\}_{\Omega}\bigr\}_{\Omega}(\mu)+\bigl
\{\phi_2,\{\phi_3,\phi_1\}_{\Omega}\bigr\}_{\Omega}(\mu)+\bigl\{\phi_3,
\{\phi_1,\phi_2\}_{\Omega}\bigr\}_{\Omega}(\mu)=0.
$$
Finally, since $d(\phi\psi)=(d\phi)\psi+\phi(d\psi)$ for any $\phi,\psi
\in C^{\infty}(\Gg)$, the Leibniz rule for the bracket is also clear. 
\end{proof}

\begin{rem}
This is the special type of Poisson bracket we will work with from 
now on.  The linear Poisson bracket $\{\ ,\ \}_{\text {lin}}$ on $\Gg=
\Gh^*$ is clearly of this type, since it corresponds to the case when
the cocycle $\Omega$ is trivial.  Meanwhile when $\Omega$ is a 
scalar-valued cocycle, we obtain the so-called {\em affine Poisson
bracket\/} \cite{Bh}, \cite{So}.  Affine Poisson structures occur 
naturally in the study of symplectic actions of Lie groups with general
(not necessarily equivariant for the coadjoint action) moment mappings.  
The notion of affine Poisson brackets has generalization also to the 
groupoid level.  See \cite{DLSW} or \cite{Wn3}.
\end{rem} 

 Suppose we are given a Poisson bracket of our special type, $\{\ ,\ \}
_{\Omega}$, on $\Gg=\Gh^*$.  To discuss its (strict) deformation
quantization, it is useful to observe that $\{\ ,\ \}_{\Omega}$ can be
viewed as a ``central extension'' of the linear Poisson bracket on the
dual vector space of the Lie algebra $\Gh/\Gz$.  This actually follows 
from the fact that the Lie bracket $[\ ,\ ]_{\Gh+V}$ given by equation
\eqref{(2.6)} can be transferred to a Lie bracket on $\Gh/\Gz\oplus V$, 
which turns out to be a central extension of the Lie bracket on $\Gh/\Gz$.
Let us make this observation more precise.

 Consider the exact sequence of Lie algebras,
$$
0\to\Gz\overset{\iota}\to\Gh\overset{\rho}\to{\Gh}/{\Gz}\to 0
$$
such that $\iota$ and $\rho$ are the injection and the quotient map,
respectively.  Let us fix a linear map $\tau:{\Gh}/{\Gz}\to\Gh$ such 
that $\rho\tau=\operatorname{id}$.  In this case, the exactness implies
that the map
\begin{equation}\label{(2.7)}
\omega_0:(x,y)\mapsto\iota^{-1}\bigl([\tau(x),\tau(y)]-\tau([x,y]
_{{\Gh}/{\Gz}})\bigr)
\end{equation}
is well-defined from ${\Gh}/{\Gz}\times{\Gh}/{\Gz}$ into $\Gz$, and it
is actually a Lie algebra cocycle for ${\Gh}/{\Gz}$ having values in
$\Gz$.  See \cite{B}, \cite{Kn}.  Then the Lie bracket on $\Gh$ 
can be written as follows:
\begin{equation}\label{(2.8)}
[X,Y]=\tau\bigl([\rho(X),\rho(Y)]_{\Gh/\Gz}\bigr)+\iota\bigl
(\omega_0(\rho(X),\rho(Y))\bigr),\qquad X,Y\in\Gh.
\end{equation}

 Since we have been regarding the center $\Gz$ as a subalgebra of $\Gh$
such that $\Gh\cap V=\Gz$, we may ignore the map $\iota$ and view $\Gz$
and its image in $\Gh$ or $V$ as the same.  Then $\tau$ is actually the 
map that determines the vector space isomorphism between $\Gh/\Gz\oplus
V$ and $\Gh+V$.  Under this isomorphism and by using equation \eqref
{(2.8)}, the Lie bracket $[\ ,\ ]_{\Gh+V}$ of equation \eqref{(2.6)}
is transferred to a Lie bracket on $\Gh/\Gz\oplus V$ defined by:
\begin{equation}\label{(2.9)}
\bigl[(x,v),(y,w)\bigr]_{\Gh/\Gz\oplus V}=[x,y]_{\Gh/\Gz}+\omega_0(x,y)
+\Omega(x,y)=[x,y]_{\Gh/\Gz}+\omega(x,y).
\end{equation}
Here $\omega_0$ is the cocycle defined in equation \eqref{(2.7)} and we
regarded $\Omega$ as a cocycle for $\Gh/\Gz$, as assured by an earlier
remark.  For convenience, we introduced a new ($V$--valued) cocycle 
$\omega$ for $\Gh/\Gz$, as a sum of the two cocycles $\omega_0$ and 
$\Omega$.  Then it is clear that equation \eqref{(2.9)} defines a central
extension of the Lie bracket $[\ ,\ ]_{\Gh/\Gz}$, where the extension is
given by the cocycle $\omega$.  We can now define a Poisson bracket on 
$\Gg$ modeled after this central extension type Lie bracket.  

\begin{theorem}\label{2.4.4}
Let $\Gh$ be a Lie algebra with center $\Gz$ and let us fix the maps
$\rho$ and $\tau$ given above.  Let $\Gg=\Gh^*$.  Consider the vector
space $V\subseteq C^{\infty}(\Gg)$ defined above and let us give $V$ 
the trivial $U(\Gh/\Gz)$--module structure.  Suppose $\omega$ is a Lie 
algebra cocycle for $\Gh/\Gz$ having values in $V$.  Then the bracket 
operation $\{\ ,\ \}_{\omega}:C^{\infty}(\Gg)\times C^{\infty}(\Gg)
\to C^{\infty}(\Gg)$ defined by
$$
\{\phi,\psi\}_{\omega}(\mu)=\bigl\langle\tau([\dot{d\phi(\mu)},\dot
{d\psi(\mu)}]_{\Gh/\Gz}),\mu\big\rangle+\omega\bigl(\dot{d\phi(\mu)},
\dot{d\psi(\mu)};\mu\bigr)
$$
is a Poisson bracket on $\Gg$.  Here $\dot{X}$ denotes the image of $X$
under the canonical projection $\rho$ of $\Gh$ onto $\Gh/\Gz$.
\end{theorem}

\begin{proof}
Define $\Omega:\Gh/\Gz\times\Gh/\Gz\to V$ by
\begin{equation}\label{(2.10)}
\Omega(x,y)=\omega(x,y)-\omega_0(x,y),
\end{equation}
where $\omega_0$ is the cocycle for $\Gh/\Gz$ defined in equation 
\eqref{(2.7)}.  Then the discussion in the previous paragraph implies 
that the bracket $\{\ ,\ \}_{\omega}$ is equivalent to the Poisson 
bracket $\{\ ,\ \}_{\Omega}$ given in Theorem \ref{2.4.3}.  Therefore,
it is clearly a Poisson bracket on $\Gg$.
\end{proof}

 Although the present formulation depends on the choice of the map 
$\tau$ and hence is not canonical, this Poisson bracket is, by 
construction, equivalent to the canonical Poisson bracket given in 
Theorem \ref{2.4.3}.  The relationship between them is given by 
equation \eqref{(2.10)}.  In particular, if we consider the cocycle 
$\omega_0$ of equation \eqref{(2.7)} in place of $\omega$, so that 
$\Omega$ is trivial, we obtain the linear Poisson bracket $\{\ ,\ \}
_{\text {lin}}$ on $\Gg$.  Therefore, to find a (strict) deformation 
quantization of our Poisson bracket $\{\ ,\ \}_{\Omega}$ of Theorem 
\ref{2.4.3}, we may as well try to find a (strict) deformation 
quantization of the central extension type Poisson bracket 
$\{\ ,\ \}_{\omega}$.  This change in our point of view is useful
when we work with specific examples, where the choice of coordinates
are usually apparent.

\section{Twisted group $C^*$--algebras as deformation quantizations}

 As we mentioned earlier, we expect that twisted group $C^*$--algebras
will be deformation quantizations of the Poisson brackets of ``central
extension'' type.  These are in fact the special type of Poisson 
brackets we defined in the previous section.  A more canonical 
description has been given in Theorem \ref{2.4.3}, while an equivalent, 
``central extension'' type description has been given in Theorem
\ref{2.4.4}.

 Let us from now on consider the Poisson bracket $\{\ ,\ \}_{\omega}$
on $\Gg=\Gh^*$, as defined in Theorem \ref{2.4.4}.  For convenience, 
we will fix the map $\tau:\Gh/\Gz\to\Gh$ and identify $\Gh/\Gz$ with 
its image $\tau(\Gh/\Gz)\subseteq\Gh$ under $\tau$.  To find a 
deformation quantization of $\{\ ,\ \}_{\omega}$, we will look for a 
group (2--)cocycle, $\sigma$, for the Lie group $H/Z$ of $\Gh/\Gz$, 
corresponding to the Lie algebra cocycle $\omega$.  Then we will form 
a twisted group $C^*$--algebra of $H/Z$ with $\sigma$, which we will 
show below will give us a strict deformation quantization of 
$C^{\infty}(\Gg)$ in the direction of $\{\ ,\ \}_{\omega}$.  By the 
equivalence of the Poisson brackets $\{\ ,\ \}_{\omega}$ and 
$\{\ ,\ \}_{\Omega}$, this may also be interpreted as giving a strict
deformation quantization of $C^{\infty}(\Gg)$ in the direction of 
$\{\ ,\ \}_{\Omega}$.  This result will be a generalization of the 
result by Rieffel \cite{Rf3} saying that an ordinary group 
$C^*$--algebra $C^*(H)$ provides a deformation quantization of 
$C^{\infty}(\Gh^*)$ in the direction of the linear Poisson bracket 
on $\Gh^*$.

 Recall that the cocycle $\omega$ provides a Lie bracket on the space 
$\Gh/\Gz\oplus V$.  If we restrict this Lie bracket to $\Gh/\Gz$, we 
obtain the map $[\ ,\ ]_{\omega}:\Gh/\Gz\times\Gh/\Gz\to\Gh/\Gz
\oplus V$ defined by:
\begin{equation}\label{(2.11)}
[x,y]_{\omega}=[x,y]_{\Gh/\Gz}+\omega(x,y).
\end{equation}
For the time being, to make our book keeping simpler, let us denote by
$\Gk$ and $K$ the Lie algebra $\Gh/\Gz$ and its Lie group $H/Z$.  We
now try to construct a group-like structure corresponding to $[\ ,\ ]
_{\omega}$.  From equation \eqref{(2.11)}, we expect to obtain a 
cocycle extension of the Lie group $K=H/Z$ via a certain group cocycle
corresponding to $\omega$.  As a first step, let us consider the 
following Baker--Campbell--Hausdorff series for $\Gk\oplus V$, 
ignoring the convergence problem for the moment.  Define
$$
S(X,Y)=X+Y+\frac12[X,Y]_{\Gk\oplus V}+\frac1{12}[X,[X,Y]_{\Gk\oplus V}
]_{\Gk\oplus V}+\frac1{12}[Y,[Y,X]_{\Gk\oplus V}]_{\Gk\oplus V}+\dots
$$
for $X,Y\in\Gk\oplus V$.  Let us also define $S_{\hbar}$ by $S_{\hbar}
(X,Y)=\frac1{\hbar}S(\hbar X,\hbar Y)$ for $\hbar\ne0$ in 
$\mathbb{R}$.  For $\hbar=0$, we let $S_0(X,Y)=X+Y$.

\begin{lem}\label{2.4.5}
Let $\hbar\in\mathbb{R}$ be fixed and let $S_{\hbar}$ be as above.  
Then we have, at least formally (ignoring the convergence problem),
\begin{align}
&S_{\hbar}\bigl(X,S_{\hbar}(Y,Z)\bigr)=S_{\hbar}\bigl(S_{\hbar}
(X,Y),Z\bigr) \notag \\
&S_{\hbar}(X,-X)=0,\quad S_{\hbar}(X,0)=S_{\hbar}(0,X)=X  \notag
\end{align}
for $X,Y,Z\in\Gk\oplus V$.
\end{lem}

\begin{proof}
Since $(X,Y)\mapsto\frac1{\hbar}[\hbar X,\hbar Y]_{\Gk\oplus V}=\hbar
[X,Y]_{\Gk\oplus V}$ is a Lie bracket, the above property of the 
Baker--Campbell--Hausdorff series is a standard result in Lie 
algebra theory \cite{B}, \cite{V}.
\end{proof}

 Note that when the cocycle $\omega$ is trivial, the map
$$
S_{\hbar}(x,y)=x*_{\hbar}y,\qquad x,y\in\Gk
$$
is an associative multiplication defined locally in a neighborhood
of $(0,0)$, on which the series converges \cite{V}.  This becomes 
a globally well-defined group multiplication on $\Gk$ when $\Gk$ 
is an exponential solvable Lie algebra.  In general when the 
cocycle $\omega$ is nontrivial, the convergence problem of the 
series $S_{\hbar}$ is not as simple because we are allowing $V$ to
be an infinite dimensional vector space.  Unless we have more 
knowledge about the Lie algebra and the cocycle, we cannot avoid 
this rather serious problem.  But let us postpone the discussion 
of the convergence problem a while longer and consider, purely 
formally, a restriction of the series $S_{\hbar}$ to $\Gk\times
\Gk$.  Let us  write:
\begin{equation}\label{(2.12)}
S_{\hbar}(x,y)=x*_{\hbar}y+R_{\hbar}(x,y)
\end{equation}
where $x,y\in\Gk$ and $x*_{\hbar}y$ is the ``group multiplication''
on $\Gk$ defined as above.  Since $[x,y]\in\Gk$ and since $[x,y]
_{\omega}-[x,y]=\omega(x,y)$ lies in $V$ which is central, it is 
clear that $R_{\hbar}(x,y)\in V$, if it converges, and this is the
term carrying all the information on the twisting of the Lie 
bracket.  It turns out that $R_{\hbar}(\ ,\ )$ is a ``group 
cocycle'' for $(\Gk,\ast_{\hbar})$ having values in $V$.

\begin{prop}\label{2.4.6}
Let the notation be as above.  Then the map $R_{\hbar}$ is a 
``group cocycle'' for $(\Gk,\ast_{\hbar})$ having values in the 
additive abelian group $V$.  That is, it satisfies, at least 
formally, the following conditions for a normalized group cocycle:
\begin{align}
&R_{\hbar}(y,z)+R_{\hbar}(x,y*_{\hbar}z)=R_{\hbar}(x,y)
+R_{\hbar}(x*_{\hbar}y,z)  \notag \\
&R_{\hbar}(x,-x)=0,\quad R_{\hbar}(x,0)=R_{\hbar}(0,x)=0 \notag
\end{align}
for $x,y,z\in\Gk$.
\end{prop}

\begin{proof}
From Lemma \ref{2.4.5}, we know that
$$
S_{\hbar}\bigl(x,S_{\hbar}(y,z)\bigr)=S_{\hbar}\bigl(S_{\hbar}
(x,y),z\bigr),\qquad x,y,z\in\Gk.
$$
If we rewrite both sides, since $R_{\hbar}(\ ,\ )$ is central, 
we have:
\begin{align}
{\text {(LHS)}}&=S_{\hbar}\bigl(x,y*_{\hbar}z+R_{\hbar}(y,z)\bigr)
=x*_{\hbar}(y*_{\hbar}z)+R_{\hbar}(y,z)+R_{\hbar}(x,y*_{\hbar}z), 
\notag \\
{\text {(RHS)}}&=S_{\hbar}\bigl(x*_{\hbar}y+R_{\hbar}(x,y),z\bigr)
=(x*_{\hbar}y)*_{\hbar}z+R_{\hbar}(x,y)+R_{\hbar}(x*_{\hbar}y,z). 
\notag
\end{align}
So we have the following equation:
$$
R_{\hbar}(y,z)+R_{\hbar}(x,y*_{\hbar}z)=R_{\hbar}(x,y)+R_{\hbar}
(x*_{\hbar}y,z).
$$
Also from the lemma, we have:
$$
R_{\hbar}(x,-x)=0,\quad R_{\hbar}(x,0)=R_{\hbar}(0,x)=0.
$$
\end{proof}

 As we mentioned above, this proposition does not make much 
sense unless we clear up the convergence problem of the series 
$S_{\hbar}$.  Note that the ``multiplication'', $\ast_{\hbar}$, 
on $\Gk$ is only locally defined.  The definition of $R_{\hbar}
(\ ,\ )$ is even more difficult because of the fact that it takes
values in an infinite dimensional vector space $V$.  Fortunately, 
in some special cases, these convergence problems do become 
simpler.  Let us mention a few here.

 When the cocycle is known to take values in a finite dimensional 
subspace $W$ of $V$, the space $\Gk+W$ becomes a finite 
dimensional Lie algebra such that the convergence of $S_{\hbar}$
is obtained at least locally in a neighborhood of 0.   But this 
case is rather uninteresting, because the cocycle extension 
becomes just another (finite dimensional) Lie algebra.  The 
corresponding Poisson bracket obtained as in Theorem \ref{2.4.4} 
is just the linear Poisson bracket on the dual space of this 
extended Lie algebra.  In particular, when $\omega=\omega_0$, the
series $S_{\hbar}$ becomes just the Baker--Campbell--Hausdorff 
series for the Lie algebra $\Gh/\Gz+\Gz=\Gh$ and the corresponding
Poisson bracket is the linear Poisson bracket on $\Gh^*$.

 Meanwhile, when the Lie algebra $\Gk$ is nilpotent, the series 
$S_{\hbar}$ becomes a finite series and hence always converges, 
whether or not the cocycle $\omega$ takes values in an infinite 
dimensional vector space.  In particular, $x*_{\hbar}y$ and 
$R_{\hbar}(x,y)$ can be defined for any $x,y\in\Gk$.

 Despite this attention to detail which we have to make, there 
are still possibilities for generalization.  There are some 
special cases of exponential solvable Lie algebras that do not 
fall into one of these cases but whose convergence problem (at 
least locally) can still be managed.  Meanwhile in a formal 
power series setting, since the convergence problem becomes less 
crucial, the result of the above proposition is still valid 
for any Lie algebra.  But in these general settings, we no 
longer expect to obtain ``strict deformation quantizations''
as we do below.  For this, some generalized notion of 
a group cocycle, in a ``local'' sense, needs to be developed.  
This search for a correct, weaker notion of deformation 
quantizations, will be postponed as a future project.

 Since we wish to establish a ``strict'' deformation quantization
(in the sense of Definition \ref{2.2.2}) of our special type of a 
Poisson bracket, we will consider the case when the Lie algebra 
$\Gk$ is nilpotent.  So from now on, let us assume that the Lie 
algebra $\Gh$ is nilpotent.  Then $\Gk=\Gh/\Gz$ also becomes 
nilpotent.  Let us choose and fix a basis for $\Gh$ (for example, 
we can take the ``Malcev basis'' \cite{CG}) such that elements of 
$\Gz$ can be written as $z=(0,z)$ and elements of $\Gk=\Gh/\Gz$ can 
be written as $x=(x,0)$.  Recall that we have defined our $V$ as the
space $C^{\infty}(\Gg/\Gq)$, where $\Gq$ is the subspace $\Gq=\Gz
^{\bot}$ of $\Gg$.  

 By Proposition \ref{2.4.6}, we obtain a group cocyle 
$R_{\hbar}$ for the nilpotent Lie group $K_{\hbar}=\bigl(\Gh/\Gz,
\ast_{\hbar}\bigr)$.  Here and from now on, if $\Gg$ is a Lie algebra
with the corresponding (simply connected) Lie group $G$, we will 
denote by $G_{\hbar}$ for the (simply connected) Lie group 
corresponding to $\Gg$ whose Lie bracket is given by $\hbar[\ ,\ ]$.   
Since $R_{\hbar}$ is a continuous function--valued cocycle having 
values in $V=C^{\infty}(\Gg/\Gq)$, it is more convenient to 
introduce instead the following (continuous) family of ordinary 
$\mathbb{T}$--valued cocycles, $r\mapsto\sigma_{\hbar}^{r}$.  Below
and through out the rest of the paper, $e(t)$ denotes the 
function $\operatorname{exp}\bigl[(2\pi i)t\bigr]$.  Also $\bar{e}
(t)=\operatorname{exp}\bigl[(-2\pi i)t\bigr]$.  The proof of the 
following proposition is immediate from Proposition \ref{2.4.6}.

\begin{prop}\label{2.4.7}
For a fixed $r\in\Gg/\Gq$, define the map $\sigma_{\hbar}^{r}
:\Gh/\Gz\times\Gh/\Gz\to\mathbb{T}$ by
$$
\sigma_{\hbar}^{r}(x,y)=\bar{e}\bigl[R_{\hbar}(x,y;r)\bigr]
=\operatorname{exp}\bigl[(-2\pi i)R_{\hbar}(x,y;r)\bigr]
$$
where $R_{\hbar}(x,y;r)$ is the evaluation at $r$ of $R_{\hbar}
(x,y)\in C^{\infty}(\Gg/\Gq)$.  Then $\sigma_{\hbar}^{r}$ is a
smooth normalized group cocycle for $K_{\hbar}$ having values in 
$\mathbb{T}$.  That is, 
\begin{align}
&\sigma_{\hbar}^{\mu}(y,z)\sigma_{\hbar}^{\mu}(x,y*_{\hbar}z)
=\sigma_{\hbar}^{\mu}(x,y)\sigma_{\hbar}^{\mu}(x*_{\hbar}y,z) \notag \\
&\sigma_{\hbar}^{\mu}(x,0)=\sigma_{\hbar}^{\mu}(0,x)=1 \notag  
\end{align}
for $x,y,z\in K_{\hbar}=\bigl(\Gh/\Gz,\ast_{\hbar}\bigr)$.  Moreover, 
if we fix $x,y\in K_{\hbar}$, then $\mu\mapsto\sigma_{\hbar}^{\mu}
(x,y)$ is a $C^{\infty}$--function from $\Gg/\Gq$ into $\mathbb{T}$.   
\end{prop}

\begin{rem}
Since the functions in $V=C^{\infty}(\Gg/\Gq)$ are invariant under the
coadjoint action of $H$, by an observation made earlier, they are also
invariant under the coadjoint action of $K_{\hbar}$.  Thus the cocycle
condition of the proposition can also be interpreted as the condition 
for a {\em normalized $\alpha$--cocycle\/} in Definition \ref{2.3.8},  
where $\alpha$ in this case is the coadjoint action of $K_{\hbar}$.  
Although this interpretation is not directly needed in our discussion 
below, this still suggests a possibility of future generalization. 
\end{rem}

 Since $\sigma_{\hbar}:r\to\sigma_{\hbar}^r$ is a continuous field
of normalized $\mathbb{T}$--cocycles 
(Definition \ref{2.3.8}) for the Lie group $K_{\hbar}$, 
it follows that we can, as in section 1, define a twisted
convolution algebra $L^1\bigl(K_{\hbar},C_{\infty}(\Gg/\Gq)\bigr)$.  
For $f,g$ in $L^1\bigl(K_{\hbar},C_{\infty}(\Gg/\Gq)\bigr)$, we have:
\begin{equation}\label{(2.13)}
(f*_{\sigma_{\hbar}}g)(y;r)=\int_{\Gh/\Gz}f(x;r)g(x^{-1}*_{\hbar}y;r)
\sigma_{\hbar}^r(x,x^{-1}*_{\hbar}y)\,dx.
\end{equation}
Here $dx$ denotes the left Haar measure for the (nilpotent) Lie 
group $K_{\hbar}=(\Gh/\Gz,\ast_{\hbar})$, which is just a fixed 
Lebesgue measure for the underlying vector space $\Gk=\Gh/\Gz$.  We 
will show below that as $\hbar$ approaches $0$, the family of twisted 
convolution algebras $\bigl\{L^1\bigl(K_{\hbar},C_{\infty}(\Gg/\Gq)
\bigr)\bigr\}_{\hbar\in\mathbb{R}}$ provides a deformation quantization
of our Poisson bracket $\{\ ,\ \}_{\omega}$ on $\Gg$.  Since we prefer 
to find a deformation at the level of continuous functions on $\Gg$, 
we need to develop suitable machinary.

 Choose a Lebesgue measure, $dX$, on $\Gh$ and Lebesgue measures, 
$dx$ and $dz$, on $\Gh/\Gz$ and $\Gz$ respectively, such that we have: 
$dX=dxdz$.  These Lebesgue measures will be Haar measures for the 
(nilpotent) Lie groups corresponding to the Lie algebras $\Gh$, $\Gh/
\Gz$, and $\Gz$.  In particular, Haar measure for $K_{\hbar}=(\Gh/\Gz,
\ast_{\hbar})$ is a Lebesgue measure $dx$ for $\Gh/\Gz$.  Meanwhile, 
we may write the dual vector space $\Gg=\Gh^*$ as $\Gg=\Gq\oplus
(\Gg/\Gq)$, a direct product of subspaces, where $\Gq=\Gz^{\bot}$.
By elementary linear algebra, we can realize $\Gq$ and $\Gg/\Gq$ as 
dual vector spaces of $\Gk=\Gh/\Gz$ and $\Gz$, respectively.  
Therefore, we are able to choose dual (Plancherel) measures, $dq$ 
and $dr$, for $\Gq$ and $\Gg/\Gq$ such that $d\mu=dqdr$ becomes a 
Plancherel measure for $\Gg=\Gh^*$.  All these measures are 
essentially Lebesgue measures.

 We then define the Fourier transform, $\mathcal F$, between the 
spaces of Schwartz functions $S(\Gh)$ and $S(\Gg)$ by
$$
({\mathcal F}f)(\mu)=\int_{\Gh}f(X)\bar{e}\bigl[\langle X,\mu\rangle
\bigr]\,dX,\qquad f\in S(\Gh)
$$
and the inverse Fourier transform, ${\mathcal F}^{-1}$, from $S(\Gg)$
to $S(\Gh)$ by
$$
({\mathcal F}^{-1}\phi)(X)=\int_{\Gg}\phi(\mu)e\bigl[\langle X,\mu
\rangle\bigr]\,d\mu,\qquad \phi\in S(\Gg).
$$
Here again, $e(t)=\operatorname{exp}\bigl[(2\pi i)t\bigr]$ and $\bar
{e}(t)=\operatorname{exp}\bigl[(-2\pi i)t\bigr]$.  Our choice of the 
Plancherel measure means that we have ${\mathcal F}^{-1}({\mathcal F}
f)=f$ for all $f\in S(\Gh)$ and ${\mathcal F}({\mathcal F}^{-1}\phi)
=\phi$ for all $\phi\in S(\Gg)$.  This is the Fourier inversion 
theorem.  Let us also define the partial Fourier transform, ${}^
{\wedge}$, from $S(\Gh/\Gz\times\Gg/\Gq)$ to $S(\Gg)=S(\Gq\times
\Gg/\Gq)$ by
$$
f^{\wedge}(q;r)=\int_{\Gh/\Gz}f(x;r)\bar{e}\bigl[\langle x,q\rangle
\bigr]\,dx.
$$
Its inverse Fourier transform ${}^{\vee}$ is similarly defined from
$S(\Gg)$ to $S(\Gh/\Gz\times\Gg/\Gq)$ by replacing $\bar{e}$ with 
$e$ in the definition.  Again, we have the Fourier inversion theorem: 
$(f^{\wedge})^{\vee}=f$ for all $f\in S(\Gh/\Gz\times\Gg/\Gq)$ and 
$(\phi^{\vee})^{\wedge}=\phi$ for all $\phi\in S(\Gg)$.

 We are now ready to state and prove our main theorem.  We show that
given our Poisson bracket $\{\ ,\ \}_{\omega}$ on $\Gg$, its strict
deformation quantization is essentially given by a family of twisted
group $C^*$--algebras.  

\begin{theorem}\label{2.4.8}
Let $\Gh$ be a nilpotent Lie algebra.  Let the notation be as above 
and let $\omega$ be a Lie algebra cocycle for $\Gh/\Gz$ having 
values in $V=C^{\infty}(\Gg/\Gq)$.  Suppose that $\Gg=\Gh^*$ is 
equipped with a Poisson bracket of our special type, $\{\ ,\ \}_
{\omega}$, as in Theorem \ref{2.4.4}.  Then there exists a dense 
(with respect to the usual $\|\ \|_{\infty}$ norm) subspace, 
${\mathcal A}\subseteq S(\Gg)$ of $C_{\infty}(\Gg)$ such that for 
a fixed $\hbar\in\mathbb{R}$, the following operation, $\times_
{\hbar}$, defined by:
\begin{equation}\label{(2.14)}
(\phi\times_{\hbar}\psi)=(\phi^{\vee}\ast_{\sigma_{\hbar}}
\psi^{\vee})^{\wedge},\qquad\phi,\psi\in{\mathcal A}
\end{equation}
is a well-defined multiplication on $\mathcal A$.  Here $\ast_
{\sigma_{\hbar}}$ is the twisted convolution defined in equation 
\eqref{(2.13)}.  Moreover, the following properties hold:
\begin{itemize}
\item We can define a suitable involution, ${}^{*_{\hbar}}$, and a 
$C^*$--norm, $\|\ \|_{\hbar}$, on ${\mathcal A}$ such that the 
$C^*$--completion of $({\mathcal A},\times_{\hbar},{}^{*_{\hbar}})$ 
with respect to $\|\ \|_{\hbar}$ defines a $C^*$--algebra $A_{\hbar}$.
\item For $\hbar\in\mathbb{R}$, the $C^*$--algebras $A_{\hbar}$ form
a continuous field of $C^*$--algebras.
\item $\bigl({\mathcal A},\times_{\hbar},{}^{*_{\hbar}},\|\ \|_{\hbar}
\bigr)_{\hbar\in\mathbb{R}}$ is a strict deformation quantization of 
${\mathcal A}\subseteq C^{\infty}(\Gg)$ in the direction of the 
Poisson bracket $(1/{2\pi})\{\ ,\ \}_{\omega}$ on $\Gg$.  In 
particular, we have:
\begin{equation}\label{(2.15)} 
\left\|\frac{\phi\times_\hbar\psi-\psi\times_\hbar\phi}{\hbar}-
\frac{i}{2\pi}\{\phi,\psi\}_{\omega}\right\|_\hbar\to 0 
\end{equation}
as $\hbar \to 0$.
\end{itemize}
\end{theorem}

\begin{proof}
(Step 1).  The twisted convolution, equation \eqref{(2.13)}, has 
been defined between functions on $\Gh/\Gz\times\Gg/\Gq$.  To 
define a multiplication between functions on $\Gg$, we use the 
partial Fourier transform to transfer the twisted convolution to 
$S(\Gg)$.

 Although $S(\Gh/\Gz\times\Gg/\Gq)\subseteq L^1\bigl(K_{\hbar},
C_{\infty}(\Gg/\Gq)\bigr)$, it is in general not true that $S(\Gh/
\Gz\times\Gg/\Gq)$ is an algebra under the twisted convolution 
$\ast_{\sigma_{\hbar}}$, unless the cocycle is trivial.  Still, at 
least on $C_c^{\infty}(\Gh/\Gz\times\Gg/\Gq)$, the 
$C^{\infty}$--functions on $\Gh/\Gz\times\Gg/\Gq$ with 
compact support, the twisted convolution is closed.  This result 
actually corresponds to a similar result in the trivial cocycle 
case (i.\,e.\ the crossed products \cite{Pd}), and the proof is 
also done by straightforward calculation.  So for the
purpose of proving the theorem, we may take $C_c^{\infty}(\Gh/\Gz
\times\Gg/\Gq)$ as the subspace on which the twisted convolution 
is closed.  

 We will let $\mathcal A$ be the image of $C_c^{\infty}(\Gh/\Gz
\times\Gg/\Gq)$ in $S(\Gg)$ under the partial Fourier transform, 
${}^{\wedge}$.  By the inverse partial Fourier transform, ${}^
{\vee}$, the subspace $\mathcal A$ is carried back onto $C_c^
{\infty}(\Gh/\Gz\times\Gg/\Gq)$.  Therefore, it follows 
immediately that equation \eqref{(2.14)} defines a closed 
multiplication on ${\mathcal A}\subseteq S(\Gg)$.  Since
$C_c^{\infty}(\Gh/\Gz\times\Gg/\Gq)$ is dense in $S(\Gh/\Gz
\times\Gg/\Gq)\subseteq L^1\bigl(K_{\hbar},C_{\infty}(\Gg/\Gq)
\bigr)$ with respect to the $L^1$--norm, it is clear that 
$\mathcal A$ is dense in $S(\Gg)\subseteq C_{\infty}(\Gg/\Gq)$ 
with respect to the $\|\ \|_{\infty}$ norm.

 Since we have defined our deformed multiplication on $\mathcal A$
to be isomorphic to the twisted convolution on $C_c^{\infty}
(\Gh/\Gz\times\Gg/\Gq)\subseteq L^1\bigl(\Gh/\Gz,C_{\infty}(\Gg/\Gq)
\bigr)$, we may also transfer other structures on the twisted 
convolution algebra to $\mathcal A$ via partial Fourier transform.
On the twisted convolution algebra, the involution is given by the 
following formula:
\begin{equation}
f^*(x;r)=\overline{f(x^{-1};r)\sigma_{\hbar}^r(x,x^{-1})}\Delta_
{K_{\hbar}}(x^{-1}).
\end{equation}
Here $\Delta_{\Gh/\Gz}\equiv1$, since $K_{\hbar}=(\Gh/\Gz,\ast_
{\hbar})$ is a nilpotent Lie group.  There also exists a canonical
$C^*$--norm, which is dominated by the $L^1$--norm of the twisted
group algebra, such that the completion with respect to the 
$C^*$--norm gives rise to the enveloping $C^*$--algebra $C^*\bigl
(K_{\hbar},C_{\infty}(\Gg/\Gq),\sigma_{\hbar}\bigr)$.  Via partial 
Fourier transform, we transfer these structures to $\mathcal A$ to 
define its involution, ${}^{*_{\hbar}}$, and the $C^*$--norm, $\|\ 
\|_{\hbar}$.  Let us denote the $C^*$--completion of $\bigl(
{\mathcal A},\times_{\hbar},{}^{*_{\hbar}},\|\ \|_{\hbar}\bigr)$ by
$A_{\hbar}$.  This proves the first assertion of the theorem.  We 
have $A_{\hbar}\cong C^*\bigl(K_{\hbar},C_{\infty}(\Gg/\Gq),\sigma_
{\hbar}\bigr)$.  Moreover, since the group $K_{\hbar}$ is amenable 
(nilpotent), the ``amenability condition'' always holds for the 
twisted convolution algebra, that is, $C^*\bigl(K_{\hbar},C_{\infty}
(\Gg/\Gq),\sigma_{\hbar}\bigr)=C^*_r\bigl(K_{\hbar},C_{\infty}
(\Gg/\Gq),\sigma_{\hbar}\bigr)$.

(Step 2: Continuity of the field of $C^*$--algebras $\{A_{\hbar}\}
_{\hbar\in\mathbb{R}}$).  For $\hbar\ne0$, there exists a group
isomorphism between $K_{\hbar}=(\Gh/\Gz,\ast_{\hbar})$ and $K=(\Gh/\Gz
,\ast)$ given by $x\mapsto\hbar x$.  For convenience, let us use the 
same notation, ${\sigma}_{\hbar}$, for the group cocycle for $K$ 
transferred by the isomorphism from the cocycle ${\sigma}_{\hbar}$ for
$K_{\hbar}$.  Then we have  $A_{\hbar}\cong C^*\bigl(K,C_{\infty}
(\Gg/\Gq),{\sigma}_{\hbar}\bigr)$, with the cocycle ${\sigma}_{\hbar}$
now viewed as the cocycle for $K$.  Moreover, it is not difficult to 
see that $\hbar\to{\sigma}_{\hbar}$ forms a continuous field of 
cocycles.  Since each $A_{\hbar}$ is a twisted group $C^*$--algebra 
satisfying the amenability condition, we conclude by Theorem \ref{2.3.9} 
that $\{A_{\hbar}\}_{\hbar\in\mathbb{R},\hbar\ne0}$ forms a continuous
field of $C^*$--algebras.

 When $h=0$, we no longer have the isomorphism between $K_0$ and $K$ 
in general.  Therefore, $A_0$ cannot be regarded as a twisted group 
algebra of $K$.  However, note that the problem will go away when 
$\Gh/\Gz$ is abelian.  In this case, each Lie group $K_{\hbar}$, 
including $\hbar=0$, is just the additive Lie group $K_0$.  That is, 
$x*_{\hbar}y=x+y$ for every $x,y\in K_{\hbar}$.  So each 
$C^*$--algebra $A_{\hbar}$ becomes $A_{\hbar}\cong C^*\bigl(K_0,
C_{\infty}(\Gg/\Gq),{\sigma}_{\hbar}\bigr)$.  Here we have used the 
notation $K_0$ instead of $K_{\hbar}$ to emphasize the fact that it is
the same additive abelian Lie group for every $\hbar\in\mathbb{R}$.  
But note that the cocycle still depends on $\hbar$.  Since $\hbar\to
{\sigma}_{\hbar}$ can be viewed as a continuous field of cocycles for 
$K_0$ and since the amenability condition holds, it follows that 
$\{A_{\hbar}\}_{\hbar\in\mathbb{R}}$ is a continuous field of 
$C^*$--algebras.

 In general, $\Gh/\Gz$ is not abelian and this argument is no longer
valid.  In this case, we may break down the nilpotent Lie algebra 
$\Gh/\Gz$ into its center and the corresponding quotient algebra.  If
the resulting quotient algebra is not abelian, we again break it into
its center and the quotient.  Since $\Gh/\Gz$ is a finite dimensional
nilpotent Lie algebra, it is clear that this process will end with
our $\Gh/\Gz$ broken into several abelian Lie algebras.  By using the
nontrivial structural theorems by Packer and Raeburn (Theorems \ref{b1}
and \ref{b2}), we can now prove the continuity of the field of 
$C^*$--algebras $\{A_{\hbar}\}_{\hbar\in\mathbb{R}}$.

 Recall that for a given $\hbar\in\mathbb{R}$,
$$
A_{\hbar}\cong C^*\bigl(K_{\hbar},C_{\infty}(\Gg/\Gq),{\sigma}
_{\hbar}\bigr)
$$
which we will write as $A_{\hbar}=B^0_{\hbar}\times_{\alpha^0(\hbar),
\sigma^0(\hbar)}(N^0)_{\hbar}$.  That is, $B^0_{\hbar}=C_{\infty}
(\Gg/\Gq)$ is the $C^*$--algebra, $N^0=K$ is the nilpotent Lie group, 
and $\sigma^0({\hbar})={\sigma}_{\hbar}$ is the cocycle for $(N^0)_
{\hbar}=K_{\hbar}$.  For the moment, $\alpha^0(\hbar)$ is the trivial 
action.  When $N^0$ is an abelian Lie group so that $(N^0)_{\hbar}=N^0$
for all $\hbar$, we have already shown that the field of $C^*$--algebras
$\{A_{\hbar}\}_{\hbar\in\mathbb{R}}$ is continuous.  We have to prove 
the result for nonabelian $N^0$.

 Since $N^0$ is a nilpotent Lie group, it has a nontrivial center $Z^0
\subseteq N^0$ as a normal subgroup.  Denote by $N^1$ the quotient Lie
group $N^0/{Z^0}$.  It is clear that $(N^1)_{\hbar}=(N^0)_{\hbar}/{(Z^0)
_{\hbar}}=(N^0)_{\hbar}/{Z^0}$, since $Z^0$ is abelian.  Since $N^1$
is also nilpotent, we can similarly define $Z^1$ and $N^2$.  We continue
this (finite) process until we have obtained an abelian group $N^k$.
Meanwhile by Theorem \ref{b1}, the $C^*$--algebra $A_{\hbar}$ can be 
written as
$$
A_{\hbar}=B^0_{\hbar}\times_{\alpha^0(\hbar),\sigma^0(\hbar)}(N^0)_
{\hbar}=\left(B^0_{\hbar}\times_{\alpha^0(\hbar),\sigma^0(\hbar)}Z^0
\right)\times_{\alpha^1(\hbar),\sigma^1(\hbar)}(N^1)_{\hbar}
$$
which we may denote by $A_{\hbar}=B^1_{\hbar}\times_{\alpha^1(\hbar),
\sigma^1(\hbar)}(N^1)_{\hbar}$.  If we define $B^2_{\hbar},B^3_{\hbar},
\dots$ in a similar manner, since $N^k$ is assumed to be abelian, we
obtain: $A_{\hbar}=B^k_{\hbar}\times_{\alpha^k(\hbar),\sigma^k(\hbar)}
N^k$.

 Let us now apply Theorem \ref{b2}.  We will prove the continuity of 
$\{A_{\hbar}\}_{\hbar\in\mathbb{R}}$ by induction on $k$.  When $k=0$, 
note that $(N^0)_{\hbar}=N^0$ since it is abelian.  This is just the 
case we have proved earlier. As an induction hypothesis, suppose that
the result holds for all positive integer less then $k$.  Since by 
definition $B^k_{\hbar}=B^{k-1}_{\hbar}\times_{\alpha^{k-1}(\hbar),
\sigma^{k-1}(\hbar)}Z^{k-1}$, where $Z^{k-1}$ is an abelian Lie group,
it follows that the field $\{B^k_{\hbar}\}_{\hbar\in\mathbb{R}}$ is 
continuous.  Let us denote by $B$ the corresponding $C^*$--algebra of
sections.  Note that $\bigl(\alpha^k(\hbar),\sigma^k(\hbar)\bigr)$ is 
a twisted action of $N^k$ on each fibre $B^k_{\hbar}$ while the 
continuity of $\hbar\to\bigl(\alpha^k(\hbar),\sigma^k(\hbar)\bigr)$ 
is obvious from the construction, which we will regard as a continuous 
field of twisted action $(\alpha,\sigma)$ on $B$.  Therefore, we 
conclude from Theorem \ref{b2} that the $C^*$--algebra $B\times_
{\alpha,\sigma}N^k$ is the algebra of sections of a $C^*$--bundle over
$\mathbb{R}$ with fibres isomorphic to $B^k_{\hbar}\times_{\alpha^k
(\hbar),\sigma^k(\hbar)}N^k$.  This means that the field of 
$C^*$--algebras $\hbar\to A_{\hbar}=B^k_{\hbar}\times_{\alpha^k(\hbar),
\sigma^k(\hbar)}N^k$ is continuous.

(Step 3: Proof of the deformation property).  On ${\mathcal A}$, we 
will form the expression, $(\phi\times_{\hbar}\psi-\psi\times_{\hbar}
\phi)/{\hbar}$, and compare this with our Poisson bracket on $\Gg$ 
defined in Theorem \ref{2.4.4}.

 Recall that a given function $\phi\in S(\Gg)$ can be written as 
$$
\phi(\mu)=\int({\mathcal F}^{-1}\phi)(X)\bar{e}\bigl[\langle X,\mu
\rangle\bigr]\,dX
$$
by the Fourier inversion theorem.  So we have
$$
d\phi(\mu)=(-2\pi i)\int({\mathcal F}^{-1}\phi)(X)\bar{e}\bigl[\langle
X,\mu\rangle\bigr]X\,dX.
$$
Therefore the Poisson bracket $\{\ ,\ \}_{\omega}$ becomes, for $\phi,
\psi\in{\mathcal A}\subseteq S(\Gg)$,
\begin{align}
\{\phi,\psi\}_{\omega}(\mu)&=\bigl\langle[\dot{d\phi(\mu)},\dot{d\psi
(\mu)}],\mu\bigr\rangle+\omega\bigl(\dot{d\phi(\mu)},\dot{d\psi(\mu)};
\mu\bigr) \notag \\
&=(-4\pi^2)\int({\mathcal F}^{-1}\phi)(X)({\mathcal F}^{-1}\psi)(Y)
\bar{e}\bigl[\langle X+Y,\mu\rangle\bigr] \notag \\
&\qquad\qquad\qquad\bigl(\langle[\dot{X},
\dot{Y}],\mu\rangle+\omega(\dot{X},\dot{Y};\mu)\bigr)\,dXdY. \notag
\end{align}
If we write an element $\mu\in\Gg=\Gq\oplus(\Gg/\Gq)$ as $\mu=(q,r)$ 
and similarly, elements $X,Y\in\Gh=\Gh/\Gz\oplus\Gz$ as $X=(x,z)$, $Y=
(y,z')$, then we obtain:
\begin{align}
\{\phi,\psi\}_{\omega}(\mu)=(-4\pi^2)&\int({\mathcal F}^{-1}\phi)(x,z)
({\mathcal F}^{-1}\psi)(y,z')\bar{e}\bigl[\langle (x+y,z+z'),(q,r)
\rangle\bigr]  \notag \\
&\quad\quad\bigl(\langle[x,y],q\rangle+\omega(x,y;r)\bigr)\,dxdzdydz'.
\notag
\end{align}

 Meanwhile, deformed multiplication on $\mathcal A$ can be written as
follows:
\begin{align}
(\phi\times_{\hbar}\psi)(q,r)&=\int\phi^{\vee}(x,r)\psi^{\vee}(x^{-1}
*_{\hbar}y,r)\sigma_{\hbar}^r(x,x^{-1}*_{\hbar}y)\bar{e}\bigl[\langle
y,q\rangle\bigr]\,dxdy \notag \\
&=\int({\mathcal F}^{-1}\phi)(x,z)({\mathcal F}^{-1}\psi)(y,z')\bar{e}
\bigl[\langle z+z',r\rangle\bigr]  \notag \\
&\quad\quad\quad\bar{e}\bigl[R_{\hbar}(x,y;r)\bigr]
\bar{e}\bigl[\langle x*_{\hbar}y,q\rangle\bigr]\,dzdz'dxdy.  \notag
\end{align}
To prove the deformation property, we must now show that the 
expression $(\phi\times_{\hbar}\psi-\psi\times_{\hbar}\phi)/{\hbar}$ 
approaches $(i/{2\pi})\{\phi,\psi\}_{\omega}$,
in the sense of equation \eqref{(2.15)}, as $\hbar\to0$.

 Since each $A_{\hbar}$ (for $\hbar\ne0$) is isomorphic to the 
twisted group $C^*$--algebra $C^*\bigl(K,C_{\infty}(\Gg/\Gq),
{\sigma}_{\hbar}\bigr)$, the $C^*$--norm $\|\ \|_{\hbar}$ is dominated
by the ``$L^1$--norm'' on $L^1\bigl(K,C_{\infty}(\Gg/\Gq)\bigr)$, 
which is actually equivalent (via the partial Fourier transform in 
$r\in\Gg/\Gq$ variable) to the $L^1$--norm on $L^1\bigl(\Gh/\Gz\times
\Gz\bigr)=L^1(\Gh)$.  It is also clear that even for $\hbar=0$, the 
$C^*$--norm $\|\ \|_{\hbar=0}$ (which is just the sup norm on 
$C_{\infty}(\Gg)$) is dominated by the $L^1$--norm on $L^1(\Gh)$, 
similarly by the Fourier transform.  Therefore to prove equation 
\eqref{(2.15)}, it is sufficient to show the convergence with respect 
to the $L^1$--norm on $L^1(\Gh)$.  But first, let us show that we have
at least the pointwise convergence valid in $\mathcal A$, thereby 
giving us a mild justification to our situation.  Since
\begin{align}
\bar{e}&\bigl[\langle x*_{\hbar}y,q\rangle+R_{\hbar}(x,y;r)\bigr]
\notag \\
&=\bar{e}\left[\langle x+y,q\rangle+\frac{\hbar}2\langle
[x,y]_{\Gh/\Gz},q\rangle+\frac{\hbar}2\omega(x,y;r)+O(\hbar^2)\right] 
\notag \\
&=\bar{e}\bigl[\langle x+y,q\rangle\bigr]+(-2\pi i)
\frac{\hbar}2\bar{e}\bigl[\langle x+y,q\rangle\bigr]\bigl(\langle
[x,y],q\rangle+\omega(x,y;r)\bigr)+O(\hbar^2),  \notag
\end{align}
we have:
\begin{align}
\lim_{\hbar\to0}&\left(\frac{\bar{e}\bigl[\langle x*_{\hbar}y,q\rangle
+R_{\hbar}(x,y;r)\bigr]-\bar{e}\bigl[\langle y*_{\hbar}x,q\rangle
+R_{\hbar}(y,x;r)\bigr]}{\hbar}\right) \notag \\
&\qquad\qquad\qquad\qquad=(-2\pi i)\bar{e}\bigl[\langle x+y,q\rangle
\bigr]\bigl(\langle[x,y],q\rangle+\omega(x,y;r)\bigr). \notag
\end{align}
From this, the pointwise convergence follows.  That is, for $\mu=(q,r)$
in $\Gg$, we have:
$$
\left(\frac{\phi\times_{\hbar}\psi-\psi\times_{\hbar}\phi}{\hbar}
\right)(\mu)-\frac{i}{2\pi}\{\phi,\psi\}_{\omega}(\mu)\longrightarrow0
$$
as $\hbar\to0$.  In a formal power series setting (for example, at 
the QUE algebra level), this kind of proof (showing pointwise 
convergence) is usually sufficient.  

 Let us now consider the convergence problem with respect to the 
$L^1$--norm on $L^1(\Gh)$, transferred to the ${\mathcal A}\subseteq
S(\Gg)$ level via Fourier transform.  For the linear Poisson bracket
on $\Gg$, Rieffel \cite{Rf3} has shown the $L^1$--convergence in 
$S(\Gh)$, and hence in $S(\Gg)$.  The idea is to give a bound for the 
$L^1$--norm for the expression, $(\phi\times_{\hbar}\psi-\psi\times_
{\hbar}\phi)/{\hbar}-(i/{2\pi})\{\phi,\psi\}_{\omega}$.  
Then by Lebesgue's dominated convergence theorem, the result follows.
In our case however, since the cocycles $R(\ ,\ )$ and $\omega(\ ,\ )$
have values in $V=C_{\infty}(\Gg/\Gq)$ and since we allowed them to be
possibly non-polynomial functions, we do not in general expect the 
convergence to take place in $S(\Gg)$.  Actually, even our deformed 
multiplication, $\times_{\hbar}$, had to be defined on a subspace 
${\mathcal A}\subseteq S(\Gg)$.

 At least in ${\mathcal A}$, we are able to find a suitable $L^1$--bound
for the above expression, since the convergence involving the cocycle 
terms can now be controlled in a compact set such that on this compact 
set, we have a uniform convergence.  So the dominated convergence 
theorem again can be applied to assure $L^1$--convergence.  Thus our 
proof is complete.   
\end{proof}

\begin{rem}
Sometimes, we are able to find a bigger space ${\mathcal A}$ in which 
the deformed multiplication is closed and the deformation property of
equation \eqref{(2.15)} is satisfied.  Since elements of the subalgebra 
${\mathcal A}$ are in a certain sense ``smooth functions'', it is 
always desirable, for noncommutative geometry purposes \cite{Cn1,Cn2},
to have as big an ${\mathcal A}$ as possible.  
\end{rem}

\begin{rem}
When $\{\ ,\ \}_{\omega}$ is the linear Poisson bracket on $\Gg=\Gh^*$
(when $\omega=\omega_0$), then we can show that:
$$
A_{\hbar}\cong C^*\bigl(K_{\hbar},C_{\infty}(\Gg/\Gq),\sigma_{\hbar}
\bigr)\cong C^*\bigl(K_{\hbar},C^*(Z),\sigma_{\hbar}\bigr)\cong C^*
(H_{\hbar}).
$$
The twisted convolution \eqref{(2.14)} on $L^1\bigl(K_{\hbar},
C_{\infty}(\Gg/\Gq),\sigma_{\hbar}\bigr)$ is just the ordinary 
convolution on $L^1(H_{\hbar})$.  So in this case, the theorem implies
that a strict deformation quantization of the linear Poisson bracket on
the dual vector space of a nilpotent Lie algebra $\Gh$ is provided by
a family of ordinary group $C^*$--algebras $\bigl\{C^*(H_{\hbar})\bigr\}
_{\hbar\in\mathbb{R}}$.  This is the result obtained in \cite{Rf3}.  
The dense subalgebra $\mathcal A$ on which the deformation takes place
is the space of Schwartz functions.
\end{rem}

 Although we have obtained the theorem assuming $\Gh$ to be nilpotent,
much of the argument will work if $\Gh$ is at least an exponential 
solvable Lie algebra.  Some technical problems come up.  First, we have 
to use the Haar measure for the group $(\Gh/\Gz,\ast_{\hbar})$, which
would no longer coincide with the Lebesgue measure on $\Gh/\Gz$, to
define the twisted convolution product.  So some modifications to the  
definitions of dual Plancherel measure and Fourier transform are necessary.
We have already mentioned the serious problems of defining the group 
cocycles and of correctly formulating the context  of strict  deformation
quantization.  

 At least in the linear Poisson bracket case, Rieffel in \cite{Rf3} has 
studied these problems when $\Gh$ is exponential, and even more general 
cases including the Lie groups which are only locally diffeomorphic to 
vector spaces, by relaxing the conditions for the deformation 
quantization.  It will be an interesting future project to find a 
correct formulation of the definitions such that the twisted group 
$C^*$--algebras arising from non-nilpotent Lie algebras (at least 
those corresponding to exponential Lie algebras) fit into the 
framework of deformation quantization, probably by allowing some 
mild relaxation on the strictness condition.  
 
 In another direction, there is a natural next step of our theorem
to study more general types of twisted group $C^*$--algebras 
(``twisted crossed products'' by Packer and Raeburn), with nontrivial 
actions as well as nontrivial cocycles, as possible deformation 
quantizations of Poisson brackets.  At present,  we do not have a
genuine example that does not degenerate into either crossed products
(only actions are nontrivial) or twisted group $C^*$--algebras with 
only cocycles nontrivial.  Still, there are some positive indications
that these $C^*$--algebras would provide a right framework for the
setting mentioned above---involving an exponential 
solvable Lie algebra and its dual vector space.


\bibliography{ref}

\providecommand{\bysame}{\leavevmode\hbox to3em{\hrulefill}\thinspace}
\begin{thebibliography}{10}

\bibitem{BFFL}
F.~Bayen, M.~Flato, C.~Fronsdal, A.~Lichnerowicz, and D.~Sternheimer,
  \emph{Deformation theory and quantization {I}, {II}}, Ann. Phys. \textbf{110}
  (1978), 61--110, 111--151.

\bibitem{Bh}
K.~H. Bhaskara, \emph{Affine {P}oisson structures}, Proc. Indian Acad. Sci.
  Math. Sci. \textbf{100} (1990), no.~3, 189--202.

\bibitem{Bl}
E.~Blanchard, \emph{D\'eformations de {$C^*$}--alg\`ebres de {H}opf}, Bull.
  Soc. Math. France \textbf{124} (1996), 141--215 (French).

\bibitem{B}
N.~Bourbaki, \emph{Elements of {M}athematics, {L}ie {G}roups and {L}ie
  {A}lgebras}, Springer-Verlag, 1989, English translation of {\em \'El\'ements
  de Math\'ematique, Groupes et Alg\`ebres de Lie\/}.

\bibitem{BuS}
R.~C. Busby and H.~A. Smith, \emph{Representations of twisted group algebras},
  Trans. Amer. Math. Soc. \textbf{149} (1970), 503--537.

\bibitem{CE}
H.~Cartan and S.~Eilenberg, \emph{Homological {A}lgebra}, Princeton Math.
  Series, no.~19, Princeton Univ. Press, 1956.

\bibitem{CP}
V.~Chari and A.~Pressley, \emph{A {G}uide to {Q}uantum {G}roups}, Cambridge
  Univ. Press, 1994.

\bibitem{Cn1}
A.~Connes, \emph{{$C^*$}--alg\'ebres et g\'eometrie diff\'erentielle}, C. R.
  Acad. Sci. Paris \textbf{290} (1980), 599--604 (French).

\bibitem{Cn2}
\bysame, \emph{Noncommutative {G}eometry}, Academic Press, 1994.

\bibitem{CG}
L.~Corwin and F.~P. Greenleaf, \emph{Representations of {N}ilpotent {L}ie
  {G}roups and {T}heir {A}pplications. {P}art 1}, Cambridge studies in advanced
  mathematics, no.~18, Cambridge Univ. Press, 1990.

\bibitem{VD}
A.~Van Daele, \emph{Quantum deformation of the {H}eisenberg group}, Proceedings
  of the Satellite Conference of ICM-90, World Scientific, Singapore, 1991,
  pp.~314--325.

\bibitem{DLSW}
P.~Dazord, J.~H. Lu, D.~Sondaz, and A.~Weinstein, \emph{Affino{\"{\i}}des de
  {P}oisson}, C. R. Acad. Sci. Paris \textbf{312} (1991), 523--527 (French).

\bibitem{Dr}
V.~G. Drinfeld, \emph{Quantum groups}, Proceedings of the International
  Congress of Mathematicians (Berkeley) (A.~M. Gleason, ed.), American
  Mathematical Society, Providence, RI, 1986, pp.~798--820.

\bibitem{FD}
J.~M.~G. Fell and R.~Doran, \emph{Representations of ${}^*$--algebras,
  {L}ocally {C}ompact {G}roups, and {B}anach ${}^*$--algebraic {B}undles}, Pure
  and Applied Mathematics, no. 125, Academic Press, 1988.

\bibitem{BJK2}
B.~J. Kahng, \emph{Non-compact quantum groups arising from {H}eisenberg type
  {L}ie bialgebras}, in preparation.

\bibitem{BJK}
\bysame, \emph{Deformation quantization of some non-compact solvable {L}ie
  groups and their representation theory}, Ph.D. thesis, University of
  California, Berkeley, 1997.

\bibitem{Kn}
A.~W. Knapp, \emph{{L}ie {G}roups, {L}ie {A}lgebras and {C}ohomology},
  Mathematical Notes, no.~34, Princeton University Press, 1988.

\bibitem{Ld}
M.~Landstad, \emph{Quantizations arising from abelian subgroups}, Int. J. Math.
  \textbf{5} (1994), 897--936.

\bibitem{LW}
J.~H. Lu and A.~Weinstein, \emph{Poisson {L}ie groups, dressing transformations
  and {B}ruhat decompositions}, J. Diff. Geom. \textbf{31} (1990), 501--526.

\bibitem{My}
J.~E. Moyal, \emph{Quantum mechanics as a statistical theory}, Proc. Camb.
  Phil. Soc. \textbf{45} (1949), 99--124.

\bibitem{PR1}
J.~Packer and I.~Raeburn, \emph{Twisted crossed products of {$C^*$}--algebras},
  Math. Proc. Cambridge Phil. Soc. \textbf{106} (1989), 293--311.

\bibitem{PR2}
\bysame, \emph{Twisted crossed products of {$C^*$}--algebras {II}}, Math. Ann.
  \textbf{287} (1990), 595--612.

\bibitem{PR3}
\bysame, \emph{On the structure of twisted group {$C^*$}--algebras}, Trans.
  Amer. Math. Soc. \textbf{334} (1992), no.~2, 685--718.

\bibitem{Pd}
G.~K. Pedersen, \emph{{$C^*$}--algebras and {T}heir {A}utomorphism {G}roups},
  London Math. Soc. Monographs, no.~14, Academic Press, London, 1979.

\bibitem{Rf2}
M.~Rieffel, \emph{Continuous fields of {$C^*$}--algebras coming from group
  cocycles and actions}, Math. Ann. \textbf{283} (1989), 631--643.

\bibitem{Rf1}
\bysame, \emph{Deformation quantization of {H}eisenberg manifolds}, Comm. Math.
  Phys. \textbf{122} (1989), 531--562.

\bibitem{Rf3}
\bysame, \emph{Lie group convolution algebras as deformation quantizations of
  linear {P}oisson structures}, Amer. J. Math. \textbf{112} (1990), 657--685.

\bibitem{Rf5}
\bysame, \emph{Some solvable quantum groups}, Operator Algebras and Topology
  (W.~B. Arveson, A.~S. Mischenko, M.~Putinar, M.~A. Rieffel, and S.~Stratila,
  eds.), Proc. OATE2 Conf: Romania 1989, Pitman Research Notes Math., no. 270,
  Longman, 1992, pp.~146--159.

\bibitem{Rf4}
\bysame, \emph{Deformation quantization for actions of {$R^d$}}, Memoirs of the
  AMS, no. 506, American Mathematical Society, Providence, RI, 1993.

\bibitem{So}
J.~M. Souriau, \emph{{S}tructure des {S}yst\`emes {D}ynamiques}, Dunod, Paris,
  1970 (French).

\bibitem{SZ}
I.~Szymczak and S.~Zakrzewski, \emph{Quantum deformations of the {H}eisenberg
  group obtained by geometric quantization}, J. Geom. Phys. \textbf{7} (1990),
  553--569.

\bibitem{V}
V.~S. Varadarajan, \emph{{L}ie {G}roups, {L}ie {A}lgebras and {T}heir
  {R}epresentations}, Springer-Verlag, New York, 1984.

\bibitem{Vy}
J.~Vey, \emph{D\'eformation du crochet de {P}oisson sur une vari\'et\'e
  symplectique}, Comm. Math. Helv. \textbf{50} (1975), 421--454 (French).

\bibitem{vN}
J.~von Neumann, \emph{Die {E}indeutigkeit der {S}chr{\" o}dingerschen
  {O}peratoren}, Math. Ann. \textbf{104} (1931), 570--578 (German).

\bibitem{Wn1}
A.~Weinstein, \emph{The local structure of {P}oisson manifolds}, J. Diff. Geom.
  \textbf{18} (1983), 523--557.

\bibitem{Wn3}
\bysame, \emph{Affine {P}oisson structures}, Intern. J. Math. \textbf{1}
  (1990), 343--360.

\bibitem{ZM}
G.~Zeller-Meier, \emph{Produits crois\'es d'une {$C^*$}--alg\`ebre par un
  groupe d'auto-morphismes}, J. Math. Pures et Appl. \textbf{47} (1968),
  101--239 (French).

\end{thebibliography}

\bibliographystyle{amsplain}

\end{document}